\documentclass[11pt]{amsart}

\usepackage{amsmath}
\usepackage{amsfonts}
\usepackage{amssymb}
\usepackage{amsthm}
\usepackage[hidelinks]{hyperref}
\newcommand{\NN}{\mathbb{N}}
\newcommand{\FF}{\mathbb{F}}
\newcommand{\CC}{\mathbb{C}}
\newcommand{\TT}{\mathbb{T}}

\newcommand{\pitilde}{\widetilde{\pi}}

\newtheorem{theorem}{Theorem}[section]
\newtheorem{proposition}[theorem]{Proposition}
\newtheorem{lemma}[theorem]{Lemma}
\newtheorem{corollary}[theorem]{Corollary}
\theoremstyle{remark}
\newtheorem{remark}[theorem]{Remark}
\newtheorem{example}[theorem]{Example}
\theoremstyle{definition}
\newtheorem{definition}[theorem]{Definition}
\numberwithin{equation}{section}

\title{Explicit Formulae for $L$-values in Positive Characteristic}
\author{Rudolph Bronson Perkins}
\email{perkins@math.univ-lyon1.edu}
\address{Institut Camille Jordan, Université Claude Bernard Lyon 1, 43 boulevard du 11 novembre 1918, 69622 Villeurbanne cedex, France}
\keywords{Positive characteristic, Pellarin's $L$-series, Bernoulli-Carlitz Numbers, Anderson-Thakur Function, Anderson generating functions}
\subjclass{MSC 11M38}

\begin{document}

\begin{abstract}
We focus on the generating series for the rational special values of Pellarin's $L$-series in $1 \leq s \leq 2(q-1)$ indeterminates, and using interpolation polynomials we prove a closed form formula relating this generating series to the Carlitz exponential, the Anderson-Thakur function, and the Anderson generating functions for the Carlitz module. We draw several corollaries, including explicit formulae and recursive relations for Pellarin's $L$-series in the same range of $s$, and divisibility results on the numerators of the Bernoulli-Carlitz numbers by monic irreducibles of degrees one and two.
\end{abstract}

\maketitle

\section{Introduction and Notation}
\subsection{The setting}
Let $\NN$ denote the set of non-negative integers. Let $\FF_q$ be the finite field of $q$ elements and characteristic $p$. Our interests lie in the arithmetic of the ring $A := \FF_q[\theta]$ of polynomials in the indeterminate $\theta$ over the finite field $\FF_q$. We let $K:= \FF_q(\theta)$ be the fraction field of $A$ and $K_\infty$ the completion of $K$ with respect to the non-archimedean absolute value $|\cdot|$ determined by the assignment $|\theta| = q$. We let $\CC_\infty$ be the completion of an algebraic closure of $K_\infty$ equipped with the unique extension of $|\cdot|$, which we also write as $|\cdot|$. Finally, we shall say a non-negative integer $n$ is \textit{written in base $q$} to mean we have written $n = \sum_{i \geq 0} n_i q^i$ with $0 \leq n_i \leq q-1$ for all $i \in \NN$, and we define the \textit{length of $n$} by $l(n) := \sum_{i \geq 0} n_i$; the base $q$ representation of $n$ is unique.

We shall work with certain rigid analytic functions defined over $\CC_\infty$ coming from the multivariate affine Tate algebra. Let $s$ be a positive integer, and define the \textit{closed unit polydisc} $B_1^s$ in $\CC_\infty^s$ as the $s$-fold product of the closed unit disc $B_1 := \{ z \in \CC_\infty : |z| \leq 1\}$.
Let $t_1\dots,t_s$ be indeterminates over $\CC_\infty$. Define the \textit{$s$-variable Tate algebra} $\TT_s$ to be the completion of the polynomial ring $\CC_\infty[t_1,\dots,t_s]$ with respect to the \textit{Gauss norm}, defined for all $f \in \CC_\infty[t_1,\dots,t_s]$ by
\[||f|| := \sup_{(z_i) \in B_1^s} |f(z_1,\dots,z_s)|,\]
where $f(z_1,\dots,z_s)$ denotes evaluation of the respective indeterminates $t_1,\dots,t_s$ in $f$ at the elements $z_1,\dots,z_s \in \CC_\infty$. More explicitly, $\TT_s$ is the $\CC_\infty$-sub-algebra
 of the formal power series ring $\CC_\infty[[t_1,\dots,t_s]]$ consisting of those power series that converge on the closed unit polydisc $B_1^s$, equipped with the Gauss norm with the same definition for these power series as for polynomials above. It is readily seen that the Gauss norm is actually a non-archimedean absolute value on $\TT_s$, and in particular $||\cdot||$ satisfies $||fg|| = ||f||\cdot||g||$ for all $f,g \in \TT_s$. When working in $\TT_1$, we will drop the subscript $1$ from our notations, so elements $f \in \TT$ are power series in $\CC_\infty[[t]]$ converging on $B_1$.
 
Finally, let $\chi_1,\dots,\chi_s$ be the $\FF_q[t_1,\dots,t_s] \subseteq \TT_s$-valued $\FF_q$-algebra morphisms on $A$ determined by $\chi_i(\theta) := t_i$, for $i = 1,\dots,s$. When in the single variable situation, we write $\chi_t(a)$ for the image of $a$ under the $\FF_q$-algebra map $A \rightarrow \FF_q[t] \subseteq \TT$ determined by $\theta \mapsto t$. The reader should be aware that everything to follow depends on our choice of $\theta$ as $\FF_q$-algebra generator of $A$.

\subsection{A brief history of Pellarin's series}
Let $A_+$ be the set of monic polynomials in $A$.

In \cite{Pe-11}, for positive integers $s$ and $n$, F. Pellarin introduced the following elements of $\TT_s$,
\begin{equation}\label{deflseries}
L(\chi_1\cdots\chi_s,n) := \sum_{a \in A_+} \frac{\chi_1(a)\cdots\chi_s(a)}{a^n} = \prod_{\substack{v \in A_+ \\ v \text{ irreducible}}} \left( 1 - \frac{\chi_1(v)\cdots\chi_s(v)}{v^n} \right)^{-1}.
\end{equation}
By \cite[Cor. 34]{AP-13}, the series in \eqref{deflseries} represent entire functions of the variables  $t_1,\dots, t_s$. We shall refer to these series as \textit{the special values of Pellarin's $L$-series}, or just \textit{Pellarin's special values} for short. 

Using his theory of deformations of vectorial modular forms, Pellarin proved \cite[Theorem 1]{Pe-11} the following identity in $\TT$,
\begin{equation}\label{Pellarinsformula}
L(\chi_t,1) = \prod_{i \geq 1} \left( 1 - \frac{\theta}{\theta^{q^i}} \right)^{-1} \prod_{j \geq 1}\left(1 - \frac{t}{\theta^{q^j}}\right).
\end{equation}
Further, for all positive integers $n \equiv 1 \mod (q-1)$, he showed \cite[Theorem 2]{Pe-11} that 
\[ L(\chi_t,n) \prod_{i \geq 1} \left( 1 - \frac{\theta}{\theta^{q^i}} \right)^{n} \prod_{j \geq 0}\left(1 - \frac{t}{\theta^{q^j}}\right)^{-1}\] 
is a rational function in the intersection $K(t)\cap\TT$, but he did not compute this rational function explicitly.

Since these initial discoveries B. Angl\`es and F. Pellarin \cite{AP-13} have given an explicit, symmetric (in $t_1,\dots,t_s$) polynomial $\Delta_{n,s} \in A[t_1,\dots,t_s]$ and shown that whenever $k$ and $s$ are positive integers such that $k \equiv s \mod (q-1)$,
\[W_{n,s} := \Delta_{n,s} L(\chi_1\dots\chi_s,n) \prod_{i \geq 1} \left( 1 - \frac{\theta}{\theta^{q^i}} \right)^{n} \prod_{l = 1}^s \prod_{j \geq 0}\left(1 - \frac{t_l}{\theta^{q^j}}\right)^{-1}\] is a polynomial in $K[t_1,\dots,t_s]$\footnote{To follow, we shall refer to $L(\chi_1\cdots\chi_s,n)$, with $k \equiv s \mod (q-1)$, as the \textit{rational special values of Pellarin's series}.}. Angl\`es and Pellarin also give a good upper bound on the total degree in $t_1,\dots,t_s$ of $W_{n,s}$, and while their result does not give the explicit polynomial in $K[t_1,\dots,t_s]$ which equals $W_{n,s}$, it holds with no bound on the number of indeterminates $s$. Their point of view stresses the importance of understanding the specializations of the variables $t_1,\dots,t_s$ in $W_{n,s}$ at roots of unity in the algebraic closure in $\CC_\infty$ of $\FF_q$ and the connection of these specializations with Goss-Dirichlet $L$-values and Thakur's Gauss sums, as defined in \cite{AP-13}. When $n = 1$, the polynomial $\Delta_{n,s}$ may be taken to equal $1$, and Angl\`es and Pellarin relate the  conjectural non-vanishing of the specialization of $W_{1,s}$ at $t_i = \theta$, for $i = 1,\dots,s$, to new divisibility information for the numerators of the Bernoulli-Carlitz numbers\footnote{See Section \ref{BCsection} for their definition.}.

In \cite[Section 4.2.1]{Perkdiss} the explicit rational function in $K(t_1,\dots,t_s)\cap\TT_s$ that $W_{n,s}$ 
equals was determined under the assumption that $1 \leq s \leq q$ and $n \equiv s \mod (q-1)$, with applications to divisibility of the Bernoulli-Carlitz numbers by degree one monic irreducibles. Note that the methods in \cite{Perkdiss} depend on \textit{Pellarin's formula} \eqref{Pellarinsformula}. 

The most recent development in the life of Pellarin's special values is the realization of Angl\`es, Pellarin and Tavares-Ribeiro in \cite{APTR} that Pellarin's series above are an instance of special values of $L$-series associated to Drinfeld modules defined over Tate algebras. This allows for the development of a good formalism for working with such series, and they set up a framework for studying the arithmetic interpretation of these series, with applications to L. Taelman's unit and class modules \cite{TaDUT}, for example. 

\subsection{Anderson generating functions} \label{AGFdef}
The results of this present work improve on those of \cite{Perkdiss}, and in order to state them, we introduce some further notation.

For non-negative integers $j$, let $D_j$ be the product of all monic polynomials in $A$ of degree $j$. The \textit{Carlitz exponential function} $e_C$ is then defined as an element of $K[[z]]$ by the formula
\[e_C(z) := \sum_{j \geq 0} \frac{z^{q^j}}{D_j}.\] It is a well-known result of Carlitz (see \cite{Car35} and \cite[Chapter 3]{Gbook}) that $e_C$ is an entire in $z$, surjective $A$-module endomorphism of $\CC_\infty$, taking the $A$-algebra action of multiplication to the Carlitz module action (see also Example \ref{carlitzmodule} below), as contained in the functional equation
\[e_C(\theta z) = e_C(z)^q + \theta e_C(z),\]
and that there is an explicit element $\pitilde$, algebraic over $K_\infty$ and unique up to a choice of root $\iota \in \CC_\infty$, (which we are making right now by naming this element) of the equation $X^{q-1} + \theta \in A[X]$, such that the kernel in $\CC_\infty$ of $e_C$ is given by $\pitilde A$. 

As a means describing the kernels of the exponential functions associated to the tensor powers of the Carlitz module, Anderson and Thakur in \cite{AT-90} introduced the function $\omega$ in $\TT$ given by
\[\omega(t) := \iota \prod_{j \geq 0} \left( 1 - \frac{t}{\theta^{q^j}} \right)^{-1},\]
and it follows from \cite[Prop. 2.5.5]{AT-90} that $\lim_{t \rightarrow \theta}(\theta - t) \omega(t) = \pitilde$. Thus one now has the explicit formula
\[\pitilde = \iota\theta\prod_{ j \geq 1} \left( 1 - \frac{\theta}{\theta^{q^j}} \right)^{-1}.\] 

From the other direction, given $\pitilde$, one may construct $\omega$, with an eye toward more general Drinfeld modules, by considering  \textit{Anderson generating functions for the Carlitz module} (AGF). Each AGF is defined, for fixed $z \in \CC_\infty$, as an element of $\TT$ by
\[f_C(z;t) := \sum_{j \geq 0} e_C\left( \frac{z}{\theta^{j+1}} \right) t^j.\]
One then has $f_C(\pitilde;t) = \omega(t)$ in $\TT$. See \cite{EGP13} for this result and a nice account of Anderson generating functions for more general Drinfeld modules defined over $A$, and note that the terminology employed here differs slightly from theirs.

In going from the generating series relation to the explicit formulae for Pellarin's series presented in the next section, the reader will benefit from the following alternate expansion for $f_C(z;t)$ due to Pellarin; see \cite[Section 4.2]{Pe-08} or \cite[Prop. 3.2]{EGP13} for a proof. For all $z \in \CC_\infty$ the following identity holds in $\TT$,
\begin{equation}\label{AGFalt}
f_C(z;t) = \sum_{j \geq 0} \frac{z^{q^j}}{D_j \cdot (\theta^{q^j} - t)}.
\end{equation}
It is apparent from \eqref{AGFalt} that $f_C(z;t)$ extends to a meromorphic function on $\CC_\infty$ in the variable $t$ with simple poles at $t = \theta^{q^j}$ with respective residues $-z^{q^j}$ for $j \in \NN$.

\subsection{The main result and explicit formulae}
We make the following convention for the remainder of this paper: Empty sums shall be understood to equal zero and empty products to equal one.

\subsubsection{Generating series for the rational special values of Pellarin's series}
Theorem \ref{themainresult} below is the main result of this paper, relating all of the functions defined so far in one formula. It may be viewed as a generalization of Carlitz' result,
\begin{equation}\label{carlitzeqn}
\frac{\pitilde}{e_C(\pitilde z)} = \sum_{a \in A} \frac{1}{z - a} = z^{-1} +  \sum_{\substack{n \geq q-1 \\ (q-1) | n}} \left(\sum_{a \in A_+} \frac{1}{a^n} \right) z^{n-1},
\end{equation}
viewed as giving a closed form to the generating series for the \textit{rational} (in $K$) \textit{Carlitz zeta values} $\zeta(0) := 1$ and  $\zeta(n) := \sum_{a \in A_+} a^{-n}$, with $n$ positive and $n \equiv 0 \mod (q-1)$, appearing on the right side above. Since $e_C \in K[[z]]$, dividing both sides of \eqref{carlitzeqn} by $\pitilde$ and making the change of variables $\pitilde z \rightarrow z$ in this same identity implies that $\zeta(n)/\pitilde^n \in K$ whenever $n \equiv 0 \mod (q-1)$.

The proof of the next result is obtained from an analysis of limits of products of certain interpolation polynomials for the characters $\chi_1, \dots, \chi_s$, and it will appear in Section \ref{MTpfsection} after some preliminary results. 

\begin{theorem}\label{themainresult}
For integers $1 \leq s \leq 2(q-1)$, and $z \in \CC_\infty$ such that $|z| < 1$, we have the following identity in $\TT_s$.
\begin{equation}\label{genserieseqn}
\sum_{a \in A} \frac{\chi_1(a)\cdots\chi_s(a)}{z-a} = \frac{\pitilde}{e_C(\pitilde z)}\prod_{i = 1}^s \frac{f_C(\pitilde z;t_i)}{\omega(t_i)}  - \sum_{\substack{(\alpha_i) \in \{0,1\}^s \\ \alpha_1+\cdots+\alpha_s = s - q}} \pitilde \prod_{i = 1}^s \frac{f_C(\pitilde z; t_i)^{\alpha_i}}{\omega(t_i)}.
\end{equation}
\end{theorem}

\begin{remark}
Observe that when $s = 1$, upon applying Pellarin's formula \eqref{Pellarinsformula}, we obtain an ``$A$-expansion'' for $f_C(\pitilde z;t)$ (in the sense that all functions on the right side of the identity below may be defined as sums or products indexed by elements of $A$),
\[f_C(\pitilde z;t) = \frac{\pitilde z}{\theta - t}\prod_{\substack{v \in A_+ \\ v \text{ irreducible}}}\left( 1 - \frac{\chi_t(v)}{v} \right) \prod_{a \in A \setminus \{0\}} \left(1 - \frac{z}{a} \right)\sum_{b \in A} \frac{\chi_t(b)}{z-b}.\]
This should be compared with formula (4.6) of El-Guindy and Papanikolas in \cite{EGP13}. 
\end{remark}

\subsubsection{Explicit formulae for Pellarin's $L$-series}
Using \eqref{AGFalt} and \eqref{carlitzeqn} and expanding both sides of \eqref{genserieseqn}, about $z = 0$ using geometric series and multiplication of power series, we obtain explicit formulae for the rational special values of Pellarin's series in $1 \leq s \leq 2(q-1)$ indeterminates. We introduce some further notation to simplify the statement of these formulas. 

For $s$-tuples $\alpha := (\alpha_i)$ and $\beta := (\beta_j)$ of non-negative integers, we write $\langle \alpha,\beta \rangle := \alpha_1\beta_1+\cdots+\alpha_s\beta_s$ for the usual inner product of $\alpha$ and $\beta$. We write $1_s := (1,1,\dots,1)$ for the $s$-tuple of all ones, and we let $|\alpha| := \alpha_1+\cdots+\alpha_s = \langle 1_s, \alpha \rangle$. For $\beta \in \NN^s$, we define $q^\beta$ for the $s$-tuple $(q^{\beta_1},\dots,q^{\beta_s})$, and finally $\beta \NN^s := \{(\beta_j \alpha_j) : \alpha \in \NN^s\}$. For example, $(0,0,\dots,0)\NN^s = \{(0,0,\dots,0)\}$, and $1_s \NN^s = \NN^s$.

\begin{corollary}\label{explicitwzeta}
For integers $1 \leq s \leq 2(q-1)$ and positive integers $k \equiv s \mod (q-1)$ the function $\pitilde^{-k}L(\chi_1\cdots\chi_s,k)\omega(t_1)\cdots\omega(t_s)$ is equal in $K(t_1,\dots,t_s)\cap\TT_s$ to
\begin{eqnarray*}
\sum_{\substack{\beta \in \NN^s \\ k-|q^\beta| \geq 0}} \frac{\zeta( k-|q^\beta|)}{\pitilde^{k-|q^\beta|}}\prod_{i = 1}^s\left( \frac{1}{D_{\beta_i}\cdot(\theta^{q^{\beta_i}} - t_i)} \right) - \sum_{\substack{\alpha \in \{0,1\}^s \\ |\alpha| = s-q}}\sum_{\substack{\beta \in \alpha\NN^s \\ \langle \alpha, q^\beta \rangle = k-1 }} \prod_{i = 1}^s\left( \frac{1}{D_{\beta_i}\cdot(\theta^{q^{\beta_i}} - t_i)} \right)^{\alpha_i}.
\end{eqnarray*}
\end{corollary}

\begin{example}
To assist the reader in parsing the notation above, we present a few examples.

Let $s = q$. Then when $k = 1$, the first sum is empty and the second sum above is equal to $-1$. When $k > 1$, the  second sum is empty. 

If $s = q+1$ and $k = 2$, then the first sum is empty, and from the second sum we obtain
\[\pitilde^{-2}L(\chi_1\cdots\chi_{q+1},2)\omega(t_1)\cdots\omega(t_{q+1}) = \sum_{i = 1}^{q+1} \frac{1}{\theta - t_i}.\] 

When the second sum above is non-empty, we observe that $1 \leq s - q \leq q-2$, and hence if $k-1 = \langle \alpha,q^\beta \rangle$ for some $\alpha \in \{0,1\}^s$ and some $\beta \in \alpha\NN^s$, then $\langle \alpha,q^\beta \rangle$ is the base $q$ representation of the integer $k-1$ and $l(k-1) = s-q$. One consequence is that the second sum does not appear if $l(k-1) \neq s-q$. This observation will be used to simplify the corollary below on divisibility of the numerators of the Bernoulli-Carlitz numbers by degree two monic irreducibles of $A$.
\end{example}

\section{The proof of Theorem \ref{themainresult}}
The proof of Theorem \ref{themainresult} will appear at the end of this section after some preliminary notation and results. 
\subsection{Twisting on $\TT_s$}
We introduce the following formalism in order to give natural proofs of certain results to follow.
\subsubsection{The twisted polynomial ring}
Let $\CC_\infty\{\tau\}$ denote the non-commutative ring of twisted polynomials over $\CC_\infty$ whose elements are polynomials in $\tau$ with coefficients in $\CC_\infty$ satisfying the twisted multiplication rule
\[\tau^i c = c^{q^i}\tau^i\]
for all $c \in \CC_\infty$ and $i \in \NN$.

\begin{example}[The Carlitz module]\label{carmoddef}
An important family of twisted polynomials is given by the \textit{Carlitz module} $\FF_q$-algebra homomorphism $\mathfrak{c} : A \rightarrow \CC_\infty\{\tau\}$ determined by
\[\theta \mapsto \mathfrak{c}_\theta := \tau + \theta. \]
For all $a \in A$, we will denote the image of $a$ under $\mathfrak{c}$ by $\mathfrak{c}_a$. For example, $\mathfrak{c}_{\theta^2} = \tau^2 + (\theta^q + \theta)\tau +\theta^2$. 
\end{example}

\subsubsection{Twisting}
For an $s$-tuple $\beta \in \NN^s$, let $c_\beta$, $t^\beta$ and $|\beta|$ have their usual meanings. Recall that, as a set, $\TT_s = \left\{\sum_{\beta \in \NN^s} c_\beta t^\beta : c_\beta \in \CC_\infty \text{ and } c_\beta \rightarrow 0 \text{ as } |\beta| \rightarrow \infty\right\}$. 

For $\phi \in \TT_s$ and $i \in \NN$, we define the $i$-th Frobenius twist of $\phi = \sum_{\beta \in \NN^s} c_\beta t^\beta$ by
\[\tau^i (\phi) : = \sum_{\beta \in \NN^s} c_\beta^{q^i} t^\beta.\]
One easily checks that $\tau^0, \tau^1, \dots$ are all continuous $\FF_q[t_1,\dots,t_s]$-algebra endomorphisms of $\TT$. We extend this definition linearly to elements in $\CC_\infty\{\tau\}$ to obtain a left $\CC_\infty\{\tau\}$-module structure on $\TT_s$, which is, for general twisted polynomials, continuous, $\FF_q[t_1,\dots,t_s]$-linear, but not multiplicative.

\begin{example}[Carlitz module action] \label{carlitzmodule}
Viewing $\CC_\infty \subseteq \TT$, we obtain the action of the Carlitz module on $\CC_\infty$ via the twisting action defined above. For example, for all $z \in \CC_\infty$ we have $\mathfrak{c}_\theta(z) = z^{q}+\theta z$, and $e_C(\theta z) = \mathfrak{c}_\theta(e_C(z))$.
\end{example}

\subsection{The Anderson-Thakur function}
The following result essentially first appears in \cite[Lemma 2.5.4]{AT-90};  see also \cite[Example 2.5]{EGP13}.
\begin{theorem}[Anderson-Thakur]\label{omegadifftheorem}
The function $\omega$ generates the free rank one $\FF_q[t]$-submodule of $\TT$ of solutions to the $\tau$-difference equation
\begin{equation}\label{omegadiffeq}
\tau(\phi) = (t - \theta)\phi.
\end{equation}
\end{theorem}

\begin{remark}
Since $\omega$ has a simple pole at $t = \theta$, we conclude that any element of $\TT$ which satisfies \eqref{omegadiffeq} and is regular at $t = \theta$ (i.e. whose power-series representation in $\CC_\infty[[t]]$ converges at $t = \theta$) must be the zero function. 
\end{remark}

The following corollary of Theorem \ref{omegadifftheorem} was first noticed by Pellarin in \cite[Section 4.2.2]{Pe-08} in much greater generality; see \cite[Lemma 29]{Pe-11} for another proof. In words, the next result states that $\omega$ is a simultaneous eigenvector for the commuting family $\{\mathfrak{c}_a\}_{a \in A}$ of $\FF_q[t]$-linear operators on $\TT$, from Example \ref{carmoddef}, with eigenvalues $\{\chi_t(a)\}_{a \in A}$.

\begin{corollary}[Pellarin]
For all $a \in A$, the following identity holds in $\TT$,
\[\mathfrak{c}_a(\omega) = \chi_t(a)\omega.\]
\end{corollary}

Theorem \ref{omegadifftheorem} gives us the means of naturally generating an important family of polynomials.

\begin{definition}
For all $d \geq 0$, we define $b_d \in A[t]$ by the equality
\[\tau^d(\omega) := b_d\omega.\]
Explicitly, for all $d \geq 0$ we have,
\[b_d(t) = \prod_{j = 0}^{d-1}(t - \theta^{q^j}).\]
\end{definition}

\begin{remark}
Through specializations of $t$ in the polynomials $b_i$ we obtain the coefficients (more precisely, their reciprocals) of the Carlitz exponential and logarithm functions, see \cite[Chapter 3]{Gbook}. For all $e \geq 0$, we have
\medskip

$\begin{array}{ll}
1. & D_e := b_e(\theta^{q^e}), \text{ and } \\
2. & \ell_e := \tau (b_e)(\theta) = \frac{d}{dt}b_{e+1}(\theta).
\end{array}$
\medskip

The polynomial $D_e$ is the product of all monic polynomials of degree $e$. The polynomial $\ell_e$ equals $(-1)L_e$, where $L_e$ is the least common multiple of all monic polynomials of degree $e$. See \cite[Prop. 3.1.6]{Gbook}.
\end{remark}

The following analytic identities will be useful to follow.

\begin{lemma}\label{bdlimitlemma}
The following identities\footnote{Note that there is a sign error in \cite[Cor. 4.1.9]{Perkdiss}.} hold in $\TT$,
\medskip

$\begin{array}{ll}
1. & \displaystyle{\lim_{e \rightarrow \infty} b_{e+1}(t)/\ell_e = -\pitilde/\omega(t)} , \\
2. & \displaystyle{\lim_{e \rightarrow \infty} b_{e}(t_1)\cdots b_{e}(t_q)/\ell_e = -\pitilde/\left(\omega(t_1)\cdots\omega(t_q)\right)}  \\
3. &  \displaystyle{\lim_{e \rightarrow \infty} b_e(t)/\ell_e = 0} . 
\end{array}$
\begin{proof}
We prove the first identity. The other two are similar.

For all $e \geq 1$ we have,
\begin{eqnarray*}
b_{e+1}(t)/\ell_{e} = -\theta \frac{\prod_{j = 0}^e \left(1 - \frac{t}{\theta^{q^j}}\right)}{\prod_{j = 1}^e \left( 1 - \frac{\theta}{\theta^{q^j}} \right)}.
\end{eqnarray*}
Comparing with the definitions of $\pitilde$ and $\omega$, letting $e \rightarrow \infty$ finishes the proof.
\end{proof}
\end{lemma}

\subsection{Carlitz polynomials and Interpolation results}
In order to state the next result, we introduce the \textit{normalized Carlitz polynomials}. 
\begin{definition}
For all positive integers $d$, let $A(d)$ be the $\FF_q$-subspace of $A$ consisting of all polynomials whose degree in $\theta$ is strictly less than $d$. Let $A(0) = \{0\}$, and note that for all $d \in \NN$, we have $|A(d)| = q^d$. 

Now for $d \in \NN$, let
\[E_d(z) := D_d^{-1}\prod_{a \in A(d)}(z - a) \in K[z].\]
\end{definition}

\begin{remark} \label{carpolyexplicit}
We mention quickly that for all $d \geq 0$ and $i = 0,\dots,d$, we have $\frac{d}{dt}b_{d+1}(\theta^{q^i}) = D_i \ell_{d-i}^{q^i}$. Hence, by \cite[Theorem 3.1.5]{Gbook} we have
\[E_d(z) = \sum_{i = 0}^d \frac{z^{q^i}}{\frac{d}{dt}b_{d+1}(\theta^{q^i})}.\]
\end{remark}

\begin{corollary}\label{Adinterpcor}
For all positive integers $d$, the polynomial 
\begin{equation}\label{interpeq1}
\Xi^{(d)}(z;t) = \sum_{j = 0}^{d-1} b_j(t) E_j(z)
\end{equation}
is the unique polynomial in $K[t,z]$ whose degree in $z$ is strictly less than $|A(d)| = q^d$ that interpolates the values of $\chi_t$ on $A(d)$.
\begin{proof}
We use Carlitz' well known formula, see \cite[Cor. 3.5.4]{Gbook},
\[\mathfrak{c}_a = \sum_{j \geq 0} E_j(a)\tau^j,\]
which holds in $A\{\tau\} \subseteq \CC_\infty\{\tau\}$ for all $a \in A$.

Let $d$ be a positive integer, and let $a \in A(d)$. We see from the previous corollary that the following identity holds in $\TT$,
\[\chi_t(a)\omega = \mathfrak{c}_a(\omega) = \sum_{j = 0}^{d-1} E_j(a) \tau^j (\omega) = \sum_{j = 0}^{d-1}b_j(t)E_j(a)\omega.\]
As $\omega$ is invertible in $\TT$ we see that 
\[\chi_t(a) = \sum_{j = 0}^{d-1}b_j(t)E_j(a)\]
holds in $\TT$, and hence in $A[t]$, for all $a \in A(d)$. Since the degree in $z$ of the polynomial in \eqref{interpeq1} is at most $q^{d-1}$ the result follows.
\end{proof}
\end{corollary}

\begin{remark}
It is interesting to observe that for $j = 0,\dots, d-1$, the coefficients $b_j(t)$ in \eqref{interpeq1} above do not depend on $d$. Thus one is led to analyze the properties of the formal series 
\begin{equation}\label{wagrep}
\Xi(z;t) := \sum_{j \geq 0} b_j(t)E_j(z).
\end{equation}
We have proven \cite[Theorem 3.4.4]{Perkdiss} that, for each finite place $v$ of $A$ and for fixed choices of $t$ in $A$ which make $\chi_t$ continuous on $A$ with respect to $v$, the formal series $\Xi(z;t)$ converges for $z$ in the $v$-adic completion $A_v$ of $A$ to the unique continuous extention of $\chi_t$ to $A_v$. For this reason, we refer to $\Xi(z;t)$ as the \textit{Wagner series} for the character $\chi_t$. Wagner's \cite[Theorem 4.2]{Wag71} gives another method by which to calculate the coefficients $b_e$, for $e \in \NN$.
\end{remark}

From Corollary \ref{Adinterpcor} we immediately deduce the following result. 

\begin{corollary}\label{interppolyid}
For all positive integers $d$, the following identity holds in $K[t,z]$,
\[\sum_{a \in A(d)} \chi_t(a)\frac{\ell_d E_d(z-a)}{z-a} = \sum_{j = 0}^{d-1} b_j(t)E_j(z).\]
\begin{proof}
One easily checks, using the explicit description of $E_d$ given in Remark \ref{carpolyexplicit} and the fact that $K[t]$ is an integral domain, that the left side of the identity above is the unique polynomial in $K[t,z]$ agreeing with $\chi_t$ on $A(d)$ whose degree in $z$ is strictly less than $|A(d)|$. 
\end{proof}
\end{corollary}

\begin{definition}
For all positive integers $d$, let $[d] := \theta^{q^d} - \theta$.
\end{definition}

We will be interested in the multiplication of the polynomials in the corollary just above. The following well-known (e.g., see \cite[Second proof of Theorem 3.1.5]{Gbook}) result will assist us.

\begin{lemma}\label{carlitzmult}
For all non-negative integers $d$, the following identity holds in $K[z]$,
\[E_d(z)^q = E_d(z) + [d+1]E_{d+1}(z).\]
\end{lemma}

The next result is due to Carlitz and is proven in \cite[Section 3.2]{Gbook}.
\begin{theorem}[Carlitz]\label{edlimexp}
For all $z \in \CC_\infty$ we have
\[\lim_{d \rightarrow \infty} \ell_d E_d(z) = \pitilde^{-1} e_C(\pitilde z).\]
\end{theorem}

\subsection{Products of interpolation polynomials}
For a function $f : A(d) \rightarrow \TT_s$ we write 
\[N^{(d)} (f)(z) := \sum_{a \in A(d)} f(a)\frac{\ell_d E_d(z-a)}{z-a} \in \TT_s[z].\]
Just as in the proof of Cor. \ref{interppolyid}, using the fact that $\TT_s$ is an integral domain, we see that $N^{(d)}(f)$ is the unique polynomial in $\TT_s[z]$ agreeing with $f$ on $A(d)$ whose degree in $z$ is strictly less than $|A(d)|$. 

Now we relate $N^{(d)}(\chi_1\cdots\chi_s)$ to the product of the polynomials $N^{(d)}(\chi_1), \dots, N^{(d)}(\chi_s)$ for $1 \leq s \leq 2(q-1)$ and positive integers $d$. This is an exercise in counting and is not the obstruction to extending Theorem \ref{themainresult} to larger $s$. Rather, it is the limits that arise in later sections that prohibit us from going outside of the range $1 \leq s \leq 2(q-1)$. We direct the reader to the example in Section \ref{exampleobstruction} below for more details.

To save writing, in the next result we will use the following notation: for $i = 1,\dots,s$ we define $N^{(d)}_i := N^{(d)}(\chi_i)$ and $b^{(d)}_i := b_d(t_i)$.

We recall that, for an $s$-tuple $\alpha := (\alpha_i) \in \{0,1\}^s$, and variables $x_1,\dots,x_s$, we have let $x^\alpha := x_1^{\alpha_1}\dots x_s^{\alpha_s}$ and $|\alpha| := \alpha_1+\cdots+\alpha_s$. Further, we let $1_s - \alpha := (1 - \alpha_i) \in \{0,1\}^s$.

\begin{proposition}\label{interplem}
For all $d \geq 1$ and $1 \leq s \leq 2(q-1)$ the following identity holds in $K[t_1,\dots,t_s,z] \subseteq \TT_s[z]$,
\begin{equation}\label{intprodeq}
N^{(d+1)}(\chi_1\cdots\chi_s) = \prod_{i = 1}^s N^{(d+1)}_i - [d+1]E_{d+1}\sum_{l = q}^s (E_d)^{l-q}\sum_{\substack{\alpha \in \{0,1\}^s \\ |\alpha| = l}} (N^{(d)})^{1_s - \alpha} (b^{(d)})^{\alpha}.
\end{equation}
\begin{proof}
Let $d$ be a fixed positive integer. 

When $1 \leq s \leq q-1$, we observe that the degrees of both the left side and right side of \eqref{intprodeq} above are strictly less than $|A(d)|$, and that they agree for $z \in A(d)$. Hence they are equal. 

Now assume $q \leq s \leq 2(q-1)$. For $i = 1,\dots,s$ we have $N^{(d+1)}_i = N^{(d)}_i + b^{(d)}_i E_d$. Thus,
\begin{eqnarray}
\prod_{i = 1}^s N^{(d+1)}_i &=& \prod_{i = 1}^s \left(N^{(d)}_i + b^{(d)}_i E_d\right) \\
 &=& \sum_{l = 0}^s (E_d)^l\sum_{\substack{\alpha \in \{0,1\}^s \\ |\alpha| = l}} (N^{(d)})^{1_s - \alpha} (b^{(d)})^{\alpha} \\
\label{degeq1} &=& \sum_{l = 0}^{q-1} (E_d)^l\sum_{\substack{\alpha \in \{0,1\}^s \\ |\alpha| = l}} (N^{(d)})^{1_s - \alpha} (b^{(d)})^{\alpha} \\
\label{degeq2} &+& \sum_{l = q}^s (E_d)^{l-q+1}\sum_{\substack{\alpha \in \{0,1\}^s \\ |\alpha| = l}} (N^{(d)})^{1_s - \alpha} (b^{(d)})^{\alpha} \\
\label{degeq3} &+& [d+1]E_{d+1}\sum_{l = q}^s (E_d)^{l-q}\sum_{\substack{\alpha \in \{0,1\}^s \\ |\alpha| = l}} (N^{(d)})^{1_s - \alpha} (b^{(d)})^{\alpha},
\end{eqnarray}
where in the last two lines we have used Lemma \ref{carlitzmult}.

We observe that for $i = 1,\dots,s$ the degree in $z$ of $N^{(d)}_i$ equals $q^{d-1}$ and of $E_d$ equals $q^d$. Hence the degree in $z$ for each term indexed by $l$ in equation \eqref{degeq1} is at most 
\[lq^d + (s-l)q^{d-1} = l(q^d - q^{d-1}) + s q^{d-1} \leq (q-1)(q^d + q^{d-1}) < |A(d+1)|, \]
where the first inequality is obtained by taking $l = (q-1)$ and $s = 2(q-1)$. Similarly, the degree in $z$ of each term indexed by $l$ in \eqref{degeq2} is at most
\[(l-q+1)q^d + (s-l)q^{d-1} = l(q^d - q^{d-1}) + s q^{d-1} - (q-1)q^d \leq (q-1)q^d < |A(d+1)|, \]
where the first inequality follows by taking $l = s = 2(q-1)$. Clearly, the degree in \eqref{degeq3} is at least $q^{d+1} = |A(d+1)|$, and the polynomial in this same equation vanishes on $A(d+1)$, as $E_{d+1}$ does. 

Thus the difference,
\[\prod_{i = 1}^s N^{(d+1)}_i - [d+1]E_{d+1}\sum_{l = q}^s (E_d)^{l-q}\sum_{\substack{\alpha \in \{0,1\}^s \\ |\alpha| = l}} (N^{(d)})^{1 - \alpha} (b^{(d)})^{\alpha}\]
is a polynomial of degree strictly less than $|A(d+1)|$ in $z$ that agrees with $\chi_1\dots\chi_s$ on $A(d+1)$. The result follows.
\end{proof}
\end{proposition}

\subsection{Wagner series for $\chi_t$ and the Anderson generating function}
Recall that we are referring to the series $\Xi(z;t)$, defined in \eqref{wagrep}, as the Wagner series for $\chi_t$ due to its $v$-adic properties, for $v$ a finite place of $A$. The following lemma provides the information we will require for convergence of the Wagner series for $\chi_t$ at the infinite place of $A$. This result essentially first appears in \cite{Perk12}. 

\begin{lemma}
For each $z \in \CC_\infty$, the series $\Xi(z;t)$ defined in \eqref{wagrep} represents an element in $\TT$. Further, this series converges at $t = \theta$ to $z$.
\begin{proof}
Let $z \in \CC_\infty$. To prove both claims it will suffice to show that for $t \in \CC_\infty$ such that $|t|<|\theta^q|$ the series in \eqref{wagrep} converges uniformly. 

Let $t$ be as in the paragraph above. Elementary non-archimedean estimates imply the existence of a constant $c(z)$, depending only on $z$ such that, for all $i \in \NN$ we have $|E_i(z)| \leq c(z)/|\ell_i|$. Similarly, one sees immediately, with our assumption on $t$, that $|b_i(t)| < |\theta^q||\ell_{i-1}|$, for all $i \in \NN$. Hence,
\[ \sum_{j \geq 0} \left| b_j(t) E_j(z) \right| \leq c(z)|\theta^q|\sum_{j \geq 0}|\ell_{j-1}/\ell_j| < \infty. \]
Thus, the series defining $\Xi(z;t)$ is the uniform limit of a sequence of polynomials in $\CC_\infty[t]$. Further this series converges at $t = \theta$, and one sees immediately that it takes the value $z$ at this point.  
\end{proof}
\end{lemma}

The next result, wherein the series $\Xi(z;t)$ and the Anderson generating functions $f_C(z;t)$ are related, was pointed out to the author by F. Pellarin, and it should also be compared with (4.6) of \cite{EGP13}. A proof appears in \cite[Theorem 1.6.6]{Perkdiss} that follows Pellarin's line of thought, using his formula \eqref{Pellarinsformula} and the formalism in \cite[Section 4]{Pe-11}. Now we give a different proof using the techniques of \cite{EGP13}. We take this route here because it allows us to remove the logical dependence of Theorem \ref{themainresult} above on Theorem 1 of \cite{Pe-11}. Thus we will obtain both Theorems 1 and 2 of \cite{Pe-11} as a consequence of Theorem \ref{themainresult} above.

\begin{theorem}\label{AGFisWIS}
For each $z \in \CC_\infty$, we have the following identity in $\TT$,
\[f_C(\pitilde z;t) = \Xi(t;z)\omega(t).\]
\begin{proof}
See \cite[Section 4]{EGP13} where it is proven that 
\[\tau(f_C)(z;t) = e_C(z) + (t - \theta) f_C(z;t).\]

Now recall that $\tau^j\omega(t) = b_j(t)\omega(t)$. Thus we also have,
\begin{eqnarray*}
\tau\left( \omega(t)\sum_{j \geq 0} b_j(t) E_j(z) \right) &=& \lim_{e \rightarrow \infty}\tau\left( \omega(t)\sum_{j = 0}^e b_j(t) E_j(z) \right)\\
&=& \omega(t)\lim_{e \rightarrow \infty}\sum_{j = 0}^e b_{j+1}(t) E_j(z)^q \\
&=& \omega(t)\lim_{e \rightarrow \infty}\left(\sum_{j = 0}^e b_{j+1}(t)E_j(z) + \sum_{j = 1}^{e+1} [j]b_{j}(t)E_{j}(z)\right) \\
&=& \omega(t)\lim_{e \rightarrow \infty}\left([e+1]b_{e+1}(t)E_{e+1}(z) + (t-\theta)\sum_{j = 0}^e b_{j}(t)E_j(z)\right) \\
&=& e_C(\pitilde z) + (t - \theta)\omega(t)\sum_{j \geq 0} b_j(t) E_j(z),
\end{eqnarray*}
by the continuity of $\tau$ and Lemmas \ref{bdlimitlemma} and \ref{carlitzmult} and Theorem \ref{edlimexp}. 
Thus the difference 
\[f_C(\pitilde z;t) - \omega(t)\sum_{j \geq 0} b_j(t) E_j(z)\] 
is is an element of $\TT$, regular at $t = \theta$ by \eqref{AGFalt} and satisfying the $\tau$ difference equation $\tau(\phi) = (t - \theta)\phi$. From this and the remark following Theorem \ref{omegadifftheorem} we deduce the result. 
\end{proof}
\end{theorem}

\subsection{Proof of Theorem \ref{themainresult}}\label{MTpfsection}
We now have all we need to prove Theorem \ref{themainresult}. 

Let $s$ be as in the statement of Theorem \ref{themainresult}, and recall from the definition that 
\[N^{d}(\chi_1\cdots\chi_s)(z) := \sum_{a \in A(d)}\chi_1(a)\cdots\chi_s(a)\frac{\ell_d E_d(z-a)}{z-a} \in K[t_1,\dots,t_s,z].\]
Thus for all $z \in \CC_\infty$ such that $|z| < 1$, by Theorem \ref{edlimexp} we have,
\[\lim_{d \rightarrow \infty}N^{d}(\chi_1\cdots\chi_s)(z) = \pitilde^{-1} e_C(\pitilde z) \sum_{a \in A} \frac{\chi_1(a)\cdots\chi_s(a)}{z-a}.\]

We now let $d$ tend to $\infty$ on the right side of \eqref{intprodeq} as follows. The right side of \eqref{intprodeq} equals 
\begin{equation}\label{mainthmlim}
\prod_{i = 1}^s N^{(d+1)}_i + \ell_{d+1}E_{d+1}\sum_{l = q}^s (\ell_d E_d)^{l-q}\sum_{\substack{\alpha \in \{0,1\}^s \\ |\alpha| = l}} (N^{(d)})^{1 - \alpha} (b^{(d)})^{\alpha} / (\ell_d)^{|\alpha|-q+1},
\end{equation}
and for $\alpha \in \{0,1\}^s$, if $|\alpha| \neq q$, then $(b^{(d)})^{\alpha} / \ell_d^{|\alpha|+1 - q} \rightarrow 0$ as $d \rightarrow \infty$ by Lemma \ref{bdlimitlemma}. Thus, using Theorems \ref{edlimexp} and \ref{AGFisWIS} we see that  as $d \rightarrow \infty$ the polynomial in \eqref{mainthmlim} converges to
\[\prod_{i = 1}^s \frac{f_{C}(\pitilde z;t_i)}{\omega(t_i)} - \frac{e_C(\pitilde z)}{\omega(t_1)\dots\omega(t_s)} \sum_{\substack{\alpha \in \{0,1\}^s \\ |\alpha| = q}} f_C(\pitilde z; t_i)^{1 - \alpha_i}.\]
This concludes the proof of Theorem \ref{themainresult}.

\subsubsection{Example obstruction in the proof of Theorem \ref{themainresult} for $s = 2(q-1)+1$} \label{exampleobstruction}
Let $s = 2q-1$. Using the same principle as in Prop. \ref{interplem}, we compute that for $d \geq 1$, the polynomial $N^{(d+1)}(\chi_1\cdots\chi_s)$ equals
\begin{eqnarray*}
\prod_{i = 1}^s N^{(d+1)}_i - \sum_{\substack{\alpha \in \{0,1\}^s \\ |\alpha| = q-1}} (b^{(d-1)})^{1_s-\alpha} (b^{(d)})^\alpha [d][d+1]E_{d+1} - (b^{(d)})^{1_s}(1+E_d^{q-1})[d+1]E_{d+1} \\
- \sum_{k = q}^{2(q-1)} \sum_{\substack{\alpha \in \{0,1\}^s \\ |\alpha| = k}} (N^{(d)})^{1_s - \alpha} (b^{(d)})^\alpha [d+1] E_d^{k-q} E_{d+1}.
\end{eqnarray*}

Now in attempting to pass to the limit as $d \rightarrow \infty$, we are lead to understand
\[\lim_{d \rightarrow \infty} \ell_{d+1}E_{d+1} \sum_{\substack{\alpha \in \{0,1\}^s \\ |\alpha| = q-1}} \left( (b^{(d-1)})^{1_s-\alpha} (b^{(d)})^\alpha / \ell_{d-1} - (b^{(d)})^{1_s} / \ell_{d} \right).\]
It is not difficult to see that the line above converges to $\pitilde^{-1} e_C(\pitilde z) L(\chi_1\cdots\chi_{2q-1},1)$, but computing its explicit value is more difficult. Angl\`es, Pellarin, and Tavares-Ribeiro give the explicit value of $L(\chi_1\cdots\chi_{2q-1},1)$ on Page 8 of \cite{APTR}, see their $\mathbb{B}_{2q-1}$.

With the expectation that such problematic limits will arise each time we increase $s$ by a multiple of $q-1$, adding limit after limit to Lemma \ref{bdlimitlemma}, combined with the combinatorial difficulty of computing a usable expression for $N^{(d)}(\chi_1\cdots\chi_s)$ for arbitrary positive integers $d$ and $s$, we leave off such pursuits for the present time.

\section{Applications to Bernoulli-Carlitz numbers}\label{BCsection}
\subsection{Preliminary Definitions and Results}
For the remainder of this note, we assume that $q \geq 3$.

\begin{definition}
For each non-negative integer $n = \sum_{i \geq 0} n_i q^i$, written in base $q$, we define the \textit{$n$-th Carlitz factorial},
\[\Pi(n) := \prod_{i \geq 0} D_i^{n_i}.\]
\end{definition}

\begin{definition}
For non-negative integers $n,k_1,\dots,k_j$ let
\begin{equation}\label{Kmultnom} \left[ \begin{array}{c} n \\ k_1,\dots,k_j \end{array} \right] := \frac{\Pi(n)}{\Pi(k_1)\Pi(k_2)\cdots\Pi(k_j)} \in K.\end{equation}
\end{definition}

The next result follows by a simple induction from \cite{Gbook} Section 9.1.
\begin{lemma}\label{carmultnomlem}
For all non-negative integers $n,k_1,\dots,k_j$ such that $k_1+\cdots+k_j = n$ the quantity defined in \eqref{Kmultnom} is an element of $A$.
\end{lemma}

Observe that the function $z/e_C(z)$ is invariant under the change of variables $z \mapsto \lambda z$ for all $\lambda \in \FF_q \setminus \{0\}$. Thus its power series representation about $z = 0$ has non-zero coefficients only for those powers of $z$ which are divisible by $(q-1)$. 

\begin{definition}
For integers $n$ divisible by $q-1$, we define the \textit{$n$-th Bernoulli-Carlitz number} $BC(n) \in K$, always assumed to have coprime numerator and denominator, via the equality
\[\frac{z}{e_C(z)} := \sum_{\substack{n \geq 0 \\ (q-1) | n}}\frac{BC(n)}{\Pi(n)}z^n \in K[[z]].\]
\end{definition}

\begin{remark}\label{explicitformulae}
From the definition of the Bernoulli-Carlitz numbers just given and \eqref{carlitzeqn} we obtain the identity $\zeta(n) / \pitilde^{n} = BC(n) / \Pi(n)$, for all non-negative integers $n$ divisible by $q-1$. Thus, when the conditions of Cor. \ref{explicitwzeta} are satisfied, we obtain the following alternate expression in $\TT_s$ for $\pitilde^{-k} L(\chi_1\cdots\chi_s,k) \omega(t_1) \cdots \omega(t_s)$,
\begin{equation*}
\sum_{\substack{\beta \in \NN^s \\ k-|q^\beta| \geq 0}} \frac{BC( k-|q^\beta|)}{\Pi({k-|q^\beta|})}\prod_{i = 1}^s\left( \frac{1}{D_{\beta_i}\cdot(\theta^{q^{\beta_i}} - t_i)} \right) - \sum_{\substack{\alpha \in \{0,1\}^s \\ |\alpha| = s-q}}\sum_{\substack{\beta \in \alpha\NN^s \\ \langle \alpha, q^\beta \rangle = k-1 }} \prod_{i = 1}^s\left( \frac{1}{D_{\beta_i}\cdot(\theta^{q^{\beta_i}} - t_i)} \right)^{\alpha_i}.
\end{equation*}
\end{remark}

See \cite[Theorem 9.2.2]{Gbook} for Carlitz' von-Staudt theorem, describing the denominators of the Bernoulli-Carlitz numbers defined above. The following corollary of Carlitz' result is all we shall require.

\begin{lemma}\label{vonstaudtlemma}
For all positive integers $n$ that are divisible by $q-1$, the denominator of $BC(n)$ either equals $1$ or else there exists a positive integer $m$ such that it is the product of all monic irreducibles of degree $m$.
\end{lemma}

\subsection{Recursive relations}
From the definition of the Bernoulli-Carlitz numbers, Carlitz observed \cite[Section 8]{Car35} that by comparing coefficients of the powers of $z$ in the formula
\[z = e_C(z)\sum_{j = 0}^\infty \frac{BC(j)}{\Pi(j)}z^j\]
one may obtain the recurrence relations $BC(0) = 1$ and for $k \geq 1$
\[BC(n) = -\sum_{q \leq q^j \leq n+1 }\left[ \begin{array}{c} n \\ n + 1 - q^j , q^j \end{array} \right] BC(n+1 - q^j),\]
giving $BC(n)$ as a $K$-linear combination of lower Bernoulli-Carlitz numbers. 

Similarly, we quickly mention that from the identity \eqref{genserieseqn}, for each $1 \leq s \leq 2(q-1)$ we obtain 
\begin{equation*}
\frac{e_C(\pitilde z)}{\pitilde}\sum_{a \in A} \frac{\chi_1(a)\cdots\chi_s(a)}{z-a} =  \prod_{i = 1}^s \frac{f_C(\pitilde z;t_i)}{\omega(t_i)}  - e_C(\pitilde z)\sum_{\substack{(\alpha_i) \in \{0,1\}^s \\ \alpha_1+\cdots+\alpha_s = s - q}} \prod_{i = 1}^s \frac{f_C(\pitilde z; t_i)^{\alpha_i}}{\omega(t_i)}.
\end{equation*}
Using the observation that the sum on the left side of the identity above is the generating series for the rational special values of Pellarin's series and expanding both sides in power series about $z = 0$, we obtain recursive relations for the rational special values of Pellarin's series by comparing the coefficients of the powers of $z$. Since we will not need this in the sequel, we leave the explicit details to the reader.

\subsection{Results for degree one monic irreducibles}
Given positive integers $s$ and $n$, one path to obtaining the Carlitz zeta values, defined  just after \eqref{carlitzeqn}, from the special values of Pellarin's series is through evaluation of the variables $t_1,\dots,t_s$ at $\theta^{q^{j_1}},\dots,\theta^{q^{j_s}}$ for non-negative integers $j_1,\dots,j_s$ such that $n - \sum_{i =1}^s q^{j_i} > 0$, and one obtains $\zeta(n - \sum_{i =1}^s q^{j_i})$ from $L(\chi_1\cdots\chi_s,n)$ in this manner. The following remark gives an alternate route that we employ below. 

\begin{remark}\label{trivcharlemma}
Let $v$ be a monic irreducible polynomial in $A$ of degree $d$, and let $\lambda \in \CC_\infty$ be one of its roots. For $a \in A$, let $a(\lambda)$ be the image of $a$ under evaluation of $\theta$ at $\lambda$. The following identity holds in $K_\infty$,
\[\sum_{a \in A_+} \frac{a(\lambda)^{q^d - 1}}{a^n} = \left( 1 - \frac{1}{v^n} \right)\zeta(n).\]
The left side of the identity above may be obtained from $L(\chi_1\cdots\chi_{d(q-1)},n)$ through specialization at $t_{i_j} = \lambda^{q^j}$ for $i_j = j(q-1) + 1, j(q-1)+2, \dots, (j+1)(q-1)$ and $j = 0,1,\dots,d-1$.
\end{remark}

In the light of this remark and the explicit formulae for the special values of Pellarin's series above, we must understand the corresponding evaluations of $\omega$ at roots of unity. We thank the anonymous reviewer for suggesting the simple algebraic proof given in the next result. 

\begin{lemma}\label{lem:omeg}
Let $v$ be a monic irreducible polynomial in $A$ of degree $d$, and let $\lambda \in \CC_\infty$ be one of its roots. The following identity holds,
\[(-1)^d v =  \prod_{i = 1}^d \omega(\lambda^{q^i})^{q-1}.\]
\begin{proof}
One may check from the definition we have given above of $\omega$ as a product, that 
\[\omega(t)^q = (t^q - \theta)\omega(t^q).\]
We know from elementary algebra that the roots of $v$ are given by $\lambda^q, \lambda^{q^2},\dots,\lambda^{q^d} = \lambda$. Substituting $\lambda^{q^i}$ into the identity above, for $i = 1,\dots,d$, we obtain
\[\omega(\lambda^{q^i})^q = (\lambda^{q^{i+1}} - \theta)\omega(\lambda^{q^{i+1}}).\]
Thus the product of all quantities on the left side of these $d$ identities equals the product of all the quantities on the right, and we obtain the result.
\end{proof}
\end{lemma}
\begin{remark}
See Angl\`es and Pellarin's results in \cite{AP-13} for the connection between the Anderson-Thakur function and Thakur's Gauss sums. Their results give a more satisfying explanation of the identity in the lemma above.
\end{remark}

We now give another recurrence relation satisfied by the Bernoulli-Carlitz numbers. We do not know the exact relation to the recursion given in the previous subsection. The benefit of the recurrence below is that when $l(n) \geq 2(q-1)$ and $\lambda \in \FF_q$ it gives $(1 - (\theta - \lambda)^n)BC(n)$ as an $A$-linear combination of lower Bernoulli-Carlitz numbers, see Remark \ref{Alincomb} below. This, in combination with Lemma \ref{vonstaudtlemma}, will allow us to obtain a lower bound on the divisibility of the numerators of the Bernoulli-Carlitz numbers by degree one monic irreducibles of $A$ as a corollary. We retain the notation of the introduction for $s$-tuples of non-negative integers.

\begin{proposition}\label{BCrecur}
We have $BC(0) = 1$, and for any integer $n$ divisible by $q-1$ and any $\lambda \in \mathbb{F}_q$,
\begin{eqnarray*} 
(1 - (\theta - \lambda)^{n})BC(n) = \sum_{\substack{\beta \in \NN^{q-1} \\ n - |q^\beta| \geq 0}} \left[ \begin{array}{c} n \\ n - |q^\beta| , |q^\beta| \end{array} \right] (\theta - \lambda)^{n -1- |q^\beta|} BC\left(n - |q^\beta|\right) .
\end{eqnarray*}
\begin{proof}
This is an elementary algebraic calculation made using the identity in Remark \ref{explicitformulae}, Remark \ref{trivcharlemma}, and Lemma \ref{lem:omeg}.
\end{proof}
\end{proposition}

\begin{remark}
This same trick can be played to obtain similar relations for $L(\chi_t,n)$ using the explicit formula for $L(\chi_1\cdots\chi_{q},n)$ given above. The details are left to the reader. 
\end{remark}

\begin{remark}\label{Alincomb}
Let us assume that $n$ is an integer that is divisible by $q-1$ and is such that $l(n) \geq 2(q-1)$. The length of $n$ is the least number of powers of $q$ needed to represent $n$. Thus for all $\beta \in \NN^{q-1}$ we have $n > |q^\beta|$. Hence, using the divisibility by $q-1$ of $n$ we have, $n - |q^\beta| \geq q-1$. It follows from Prop. \ref{BCrecur} and Lemma \ref{carmultnomlem} that $(1 - (\theta - \lambda)^n)BC(n)$ may be given as an $A$-linear combination of lower Bernoulli-Carlitz numbers.
\end{remark}

The following notation will allow us to better state the lower bound on divisibility of the numerators of the Bernoulli-Carlitz numbers by degree one monic irreducibles in $A$, and we will use it in the next subsection. 

\begin{definition}\label{1norm}
Let $n$ and $s$ be positive integers. For $\beta \in \NN^s$, define
\[|\beta|_1 := |q^\beta| = \sum_{i = 1}^s q^{\beta_i}.\]

Let $\mathcal{M}_s(n)$ be the set of tuples $\beta \in \NN^{s(q-1)}$ such that $n \geq |\beta|_1$. We call a tuple $\mu \in \mathcal{M}_s(n)$ \textit{maximal for $|\cdot|_1$} if $|\mu|_1 \geq |\beta|_1$ for all $\beta \in \mathcal{M}_s(n)$.
\end{definition}

\begin{corollary}\label{BCdegone}
Let $n$ be a positive integer divisible by $q-1$ and such that $l(n) \geq 2(q-1)$, and let $\mu \in \mathcal{M}_1(n)$ be maximal for $|\cdot|_1$. If the denominator of $BC(n - |\mu|_1)$ is not $\theta^q - \theta$, then $(\theta^q - \theta)^{n-1-|\mu|_1}$ divides the numerator of $BC(n)$. Otherwise, $(\theta^q - \theta)^{n-2-|\mu|_1}$ divides the numerator of $BC(n)$.
\begin{proof}
Let $n$ be as above, and let $\lambda \in \FF_q$. Using Lemma \ref{BCrecur} we see that $(1 - (\theta - \lambda)^{n})BC(n)$ equals
\begin{equation*}
\sum_{\beta \in \mathcal{M}_1(n)} \left[ \begin{array}{c} n \\ n - |\beta|_1, |\beta|_1 \end{array} \right] (\theta - \lambda)^{n - 1 - |\beta|_1} BC\left(n - |\beta|_1 \right).
\end{equation*}
By our assumptions on $n$ and Remark \ref{Alincomb} all the differences $n - 1 - |\beta|_1$ appearing in the exponent of $\theta - \lambda$ above are at least $q-2$, and $n - 1 - |\mu|_1$ is minimal among them. 
Taking account of the possibility that $\theta^q - \theta$ equals the denominator of some $BC(j)$ appearing above and using that $\lambda \in \mathbb{F}_q$ was arbitrary finishes the proof.
\end{proof}
\end{corollary}

\begin{remark} By Remark \ref{Alincomb}, outside of the possibility that $q = 3$ and that the denominator of $BC(n - |\mu|_1)$ equals $\theta^3 - \theta$, the corollary above always gives non-trivial information. See the examples subsection below for improvements and other remarks.
\end{remark}

\subsection{Results for degree two monic irreducibles}
In the next result, we will focus on the identity in Remark \ref{explicitformulae} when $s = 2(q-1)$ and specialize at $t_1,\dots,t_{q-1} = \lambda$ and $t_q,\dots,t_{2(q-1)} = \lambda^q$, for $\lambda$ a root of a monic irreducible polynomial in $A$ of degree two. We remind the reader that $1_{2(q-1)} \in \{0,1\}^{2(q-1)}$ denotes the vector, all of whose entries equal 1. For convenience, we introduce the notation $1_{1,0} \in \{0,1\}^{2(q-1)}$ for the vector whose first $q-1$ entries are ones and whose remaining entries are zeros, and we let $1_{0,1} := 1_{2(q-1)} - 1_{1,0}$. Further, juxtaposition of two vectors $\alpha, \beta \in \NN^{2(q-1)}$ will mean their coordinate-wise product,
\[\alpha\beta := (\alpha_1\beta_1,\alpha_2\beta_2,\dots,\alpha_{2(q-1)}\beta_{2(q-1)}).\]
For example, $q^\beta 1_{1,0} = (q^{\beta_1},\dots,q^{\beta_{q-1}},0,\dots,0)$ and $q^\beta 1_{0,1} = (0,\dots,0,q^{\beta_q},\dots,q^{\beta_{2(q-1)}})$. 

\begin{proposition}\label{degtworecur}
Let $\lambda \in \CC_\infty$ be the root of a monic irreducible polynomial $v \in A$ of degree two. For each integer $n$ divisible by $q-1$ and such that $l(n) \geq 2(q-1)$ we have
\begin{eqnarray*} 
(v^{n}-1)BC(n) = \sum_{\beta \in \mathcal{M}_2(n)}  \frac{\left[ \begin{array}{c} n \\ n - |q^\beta| , |q^\beta 1_{1,0}| , |q^\beta 1_{0,1}| \end{array} \right] v^{n - 1} }{\prod_{i = 1}^{q-1}(\theta^{q^{\beta_i}} - \lambda) \prod_{i = q}^{2(q-1)}(\theta^{q^{\beta_i}} - \lambda^q)}BC\left(n - |q^\beta|\right) .
\end{eqnarray*}
\begin{proof}
This follows immediately Remark \ref{explicitformulae}, Remark \ref{trivcharlemma}, Lemma \ref{lem:omeg}, and the definition of the Carlitz factorials.
\end{proof}
\end{proposition}

\begin{remark} \label{degtworemark}
For $\lambda,v = (\theta - \lambda)(\theta-\lambda^q)$ and $n$ as in the proposition just above, when $l(n) \geq 3(q-1)$, we have $n > |\beta|_1$, for all $\beta \in \mathcal{M}_2(n)$. Hence, for all $\beta \in \mathcal{M}_2(n)$, both 
\[ \left[ \begin{array}{c} n \\ n - |q^\beta| , |q^\beta 1_{1,0}| , |q^\beta 1_{0,1}| \end{array} \right] \ \text{ and } \ \frac{v^{n-1}}{\prod_{i = 1}^{q-1}(\theta^{q^{\beta_i}} - \lambda) \prod_{i = q}^{2(q-1)}(\theta^{q^{\beta_i}} - \lambda^q)}\]
are elements of $\FF_q(\lambda)[\theta]$. Thus the previous result gives $(v^n-1)BC(n)$ as an integral linear combination of lower Bernoulli-Carlitz numbers, whenever $l(n) \geq 3(q-1)$.
\end{remark}

Let $\lambda$ be as in Prop. \ref{degtworecur}. In the next result, we shall need to determine the highest power of $\theta - \lambda$ (and, by symmetry, of $\theta - \lambda^q$) that can appear in the product
\[\prod_{i = 1}^{q-1}(\theta^{q^{\beta_i}} - \lambda) \prod_{i = q}^{2(q-1)}(\theta^{q^{\beta_i}} - \lambda^q)\]
as $\beta = (\beta_i)$ ranges over all tuples in $\mathcal{M}_2(n)$.

\begin{definition}\label{2norm}
A tuple $\beta \in \NN^{2(q-1)}$ shall be called \textit{ordered} if there exists an integer $e$ in the range $0 \leq e \leq 2(q-1)$ such that $\beta_1 \geq \beta_2 \geq \cdots \geq \beta_e$ with $\beta_1,\dots,\beta_e$ all even ($0$ is considered even here), and $\beta_{e+1} \leq \beta_{e+2} \leq \dots \leq \beta_{2(q-1)}$ with $\beta_{e+1},\dots,\beta_{2(q-1)}$ all odd.

Now, let $\beta \in \NN^{2(q-1)}$ be an ordered tuple. Define
\[|\beta|_{2} := \sum_{i = 1}^{\min\{e,q-1\}}q^{\beta_i} + \sum_{i = \max\{q, e+1\}}^{2(q-1)}q^{\beta_i}.\] 

Given an integer $n$, we shall call an ordered tuple $\mu \in \mathcal{M}_2(n)$ \textit{maximal for $|\cdot|_2$} if $|\mu|_{2} \geq |\beta|_{2}$ for all ordered tuples $\beta \in \mathcal{M}_2(n)$.
\end{definition}

We may now state the lower bound on divisibility of the numerators of the Bernoulli-Carlitz numbers by degree two monic irreducibles. 

\begin{corollary}\label{BCdegtwo}
Let $P_2$ be the product of all monic irreducible polynomials in $A$ of degree two.
Let $n$ be a positive integer divisible by $q-1$ and such that $l(n) \geq 3(q-1)$. Let $\mu \in \mathcal{M}_2(n)$ be an ordered tuple that is maximal for $|\cdot|_2$.

If the denominator of $BC(n - |q^{\mu}|)$ is not $P_2$, then the numerator of $BC(n)$ is divisible by $P_2^{n - 1 - |\mu|_{2}}$. Otherwise, the the numerator of $BC(n)$ is divisible by $P_2^{n - 2 - |\mu|_{2}}$.
\begin{remark}
In the appendix, for each integer $n$ as above, we construct an explicit ordered $\mu \in \mathcal{M}_2(n)$ that is maximal for $|\cdot|_2$, and it follows that $|\mu|_2 = |\mu|_1 = |q^\mu|$.
\end{remark}
\begin{proof}
Let $\lambda \in \CC_\infty$ be the root of a monic irreducible polynomial $v \in A$ of degree two. The assumptions of Prop. \ref{degtworecur} are satisfied and we have
\begin{eqnarray*} 
(v^{n}-1)BC(n) = \sum_{\substack{\beta \in \NN^{2(q-1)} \\ n - |q^\beta| \geq 0}}  \frac{\left[ \begin{array}{c} n \\ n - |q^\beta| , |q^\beta 1_{1,0}| , |q^\beta 1_{0,1}| \end{array} \right] v^{n - 1} }{\prod_{i = 1}^{q-1}(\theta^{q^{\beta_i}} - \lambda) \prod_{i = q}^{2(q-1)}(\theta^{q^{\beta_i}} - \lambda^q)} BC\left(n - |q^\beta|\right).
\end{eqnarray*}
Remark \ref{degtworemark} implies that the right side above is an $\FF_q(\lambda)[\theta]$-linear combination of Bernoulli-Carlitz numbers. We must determine the highest power of $\theta - \lambda$ appearing in the product
\[\prod_{i = 1}^{q-1}(\theta^{q^{\beta_i}} - \lambda) \prod_{i = q}^{2(q-1)}(\theta^{q^{\beta_i}} - \lambda^q),\]
as $\beta = (\beta_i)$ ranges over all tuples in $\mathcal{M}_2(n)$. From elementary algebra we see that $(\theta^{q^{\beta_i}} - \lambda) = (\theta - \lambda)^{q^{\beta_i}}$ and $(\theta^{q^{\beta_i}} - \lambda^q) = (\theta - \lambda^q)^{q^{\beta_i}}$ if $\beta_i$ is even, and that $(\theta^{q^{\beta_i}} - \lambda^q) = (\theta - \lambda)^{q^{\beta_i}}$ and $(\theta^{q^{\beta_i}} - \lambda^q) = (\theta - \lambda)^{q^{\beta_i}}$ if $\beta_i$ is odd. Thus, for each positive integer $k$, among tuples $\beta \in \mathcal{M}_2(n)$ such that $|\beta|_1 = k$ we see that the highest power of $\theta - \lambda$ appearing in the product above comes from an ordered tuple, and we see that for such ordered $\beta$ the power of $\theta - \lambda$ appearing in the product above equals $|\beta|_2$. Hence, as $\mu$ is maximal for $|\cdot|_2$, we see that the greatest power of $\theta - \lambda$ appearing in the denominator of any summand above comes from $\mu$. Taking account of the fact that the denominator of $BC(n - |q^\beta|)$ might equal $P_2$ for some $\beta \in \mathcal{M}_2(n)$ yields the result for $\theta - \lambda$ in place of $P_2$, and because $\lambda$ and $v$ were arbitrary, the proof is finished.
\end{proof}
\end{corollary}

We immediately obtain the following weaker divisibility result for degree two monic irreducibles of $A$ which has the flavor of the divisibility result for degree one monic irreducibles.

\begin{corollary}
Let $n$ be a positive integer divisible by $q-1$ and such that $l(n) \geq 3(q-1)$. Let $\mu \in \mathcal{M}_2(n)$ be maximal for $|\cdot|_1$. The polynomial $P_2^{n - 2 - |\mu|_1}$ divides the numerator of $BC(n)$.
\begin{proof}
Let $n$ and $\mu$ be as in the statement of the result. The claim follows immediately from the previous corollary and the fact that $|\mu|_1 \geq |\beta|_1 \geq |\beta|_2$ for all $\beta \in \mathcal{M}_2(n)$.
\end{proof}
\end{corollary}

\begin{corollary}
Let $n$ be a positive integer, divisible by $q-1$. There exists an infinite sequence $(\gamma_j)$ of positive integers, divisible by $q-1$ (explicitly constructed in the proof) that converges in the $p$-adic integers to $n$ such that $BC(\gamma_j)$ converges to zero in the $v$-adic completion of $A$ for any monic irreducible polynomial $v$ of degree one or two.
\begin{proof}
Let $l$ be a positive integer such that $q^l > n$. Define
\[\gamma_j := (q-1)q^j(q^{l+2} + q^{l+1} + q^{l}) + n.\]
Clearly, this converges to $n$ in the $p$-adic integers as $j \rightarrow \infty$.

In the case $v$ is of degree one, Cor. \ref{BCdegone} gives that the numerator of $BC(\gamma_j)$ is divisible at least by $v^{(q-1)q^j(q^{l+1}+q^l)+ n - 2}$, for all $j \geq 0$. Thus $BC(\gamma_j)$ converges to $0$ with respect to the absolute value arising from $v$.

In case $v$ is of degree two, Cor. \ref{BCdegtwo} gives that the numerator of $BC(\gamma_j)$ is divisible at least by $v^{(q-1)q^j(q^l) + n - 2}$ for all $j \geq 0$, giving the result. 
 \end{proof}
\end{corollary}

\subsection{Examples}
In the following examples, we let $q = 3$ and $P_2$ be the product of all monic irreducibles in $\FF_3[\theta]$ of degree 2. Both examples below deal with divisibility of the numerators of Bernoulli-Carlitz numbers by $P_2$, but similar examples and remarks can be made for divisibility by degree one monic irreducibles.

\begin{example}
Let $n=304 = 3^5 + 2 \cdot 3^3 + 2 \cdot 3 + 1$. We have that the ordered $|\cdot|_1$-maximal tuple $\beta = (1,3,3,5)$ does not equal the ordered $|\cdot|_2$-maximal tuple $\mu = (2,2,3,5)$, and $n - |\mu|_2 = 16$. Observe that, $P_2$ equals the denominator of $BC(16)$, by Goss' table \cite[pp. 354-358]{Gbook}. Thus Cor. \ref{BCdegtwo} gives a lower bound of $n - |\mu|_2 - 2 = 14$, and we compute, using Sage, that ${14}$ is the exact power of $P_2$ dividing the numerator of $BC(304)$.
\end{example}

\begin{example}
Let $n = 646 = 2\cdot 3^5 + 3^4 + 2 \cdot 3^3 + 2 \cdot 3^2 + 2 \cdot 3 + 1$. Then the $|\cdot|_2$-maximal tuple is $\mu = (4,2,5,5)$, and $n - |\mu|_2 = 70$. Goss' table \cite[pp. 354-358]{Gbook} shows that $BC(70) \in A$, and hence Cor. \ref{BCdegtwo} gives a lower bound of $69$ as the power of $P_2$ dividing the numerator of $BC(646)$, while Sage shows that the exact power of $P_2$ dividing $BC(646)$ is $74$. We observe that Cor. \ref{BCdegtwo} also applies to $n' = 70$ and gives a lower bound of $5$ as the power of $P_2$ dividing the numerator of $BC(70)$. Of course, $5+69 = 74$, and this suggests that the lower bound of Cor. \ref{BCdegtwo} can be sharpened by iteratively taking account of the power of $P_2$ in the numerator of $BC(n - |\mu|_2)$, with $n$ and $\mu$ as in the statement of the corollary. In this direction, one must describe the second greatest power of $\theta - \lambda$ appearing in the product
\[\prod_{i = 1}^{q-1}(\theta^{q^{\beta_i}} - \lambda) \prod_{i = q}^{2(q-1)}(\theta^{q^{\beta_i}} - \lambda^q),\]
as $\beta$ ranges over $\mathcal{M}_2(n)$, and so on. For the sake of space, we do not pursue this here.
\end{example}

\begin{example}
We have confirmed the following statement experimentally for the first $400$ non-zero Bernoulli-Carlitz numbers for monic irreducibles of degrees three and four in $\FF_3[\theta]$ using Sage mathematical software, and we wonder if it holds for all positive integers $d$ and any integer $q \geq 3$.
\medskip

\textit{Let $d$ be a positive integer, and let $P_d$ be the product of all monic irreducible polynomials in $A$ of degree $d$. Let $n$ be a positive integer divisible by $q-1$ and such that $l(n) \geq (d+1)(q-1)$. Let $\mu \in \mathcal{M}_d(n)$ be maximal for $|\cdot|_1$.  The polynomial $P_d^{n - 2 - |\mu|_1}$ divides the numerator of $BC(n)$.}
\end{example}

\appendix
\section{An explicit maximal tuple for $|\cdot|_2$}
We direct the reader to Definitions \ref{1norm}, \ref{2norm} and the introduction for the notation used to follow. 
\begin{lemma}
Let $n$ be a positive integer divisible by $q-1$ and such that $l(n) \geq 3(q-1)$. If $\mu \in \mathcal{M}_2(n)$ is maximal for $|\cdot|_1$, then $\mu_j \neq 0$ for $j = 1,2,\dots,2(q-1)$.
\begin{proof}
It is easy to see that one obtains an unordered tuple $\gamma \in \mathcal{M}(n)$, maximal for $|\cdot|_1$, from the base $q$ expansion of $n = \sum_{j = 0}^d n_i q^i$ (we assume $n_d \neq 0$) by letting 
\begin{eqnarray*}
\gamma_1 = \gamma_2 = \cdots = \gamma_{n_d} &=& d \\
\gamma_{n_d+1} = \gamma_{n_d+2} = \cdots = \gamma_{n_d + n_{d-1}} &=& d-1 \\
&\vdots&
\end{eqnarray*}
and so on, until one has defined $\gamma_1,\gamma_2,\dots,\gamma_{2(q-1)}$. Since, by construction, there is no carry over in the sum $|\gamma|_1 + (n - |\gamma|_1)$, we have $l(n) = l(|\gamma|_1) + l(n - |\gamma|_1)$. We deduce that $l(n - |\gamma|_1) \geq q-1$, and hence no entry of $\gamma$ can be zero. Now all other $\mu \in \mathcal{M}_2(n)$ that are maximal for $|\cdot|_1$ are obtained from $\gamma$ by permuting the entries. 
\end{proof}
\end{lemma}

Let $n$ be as in the previous lemma. Taking the explicit tuple $\gamma \in \mathcal{M}_2(n)$ constructed in the previous lemma, and ordering its entries as in Definition \ref{2norm}, we obtain an ordered tuple $\beta$ that is maximal for $|\cdot|_1$ (in the sense of Definition \ref{1norm}). In the next lemma we use $\beta$ to explicitly give an ordered tuple $\mu \in \mathcal{M}_2(n)$ that is maximal for $|\cdot|_2$.

\begin{proposition}
Let $n$ be a positive integer divisible by $q-1$ and such that $l(n) \geq 3(q-1)$, let $\beta \in \mathcal{M}_2(n)$ be an ordered tuple which is maximal for $|\cdot|_1$, and let $e$ be the number even entries. If $e = q-1$, then $|\beta|_1 = |\beta|_2$, and $\beta$ is maximal for $|\cdot|_2$. Otherwise, suppose that $0 \leq e < q-1$ (respectively, that $q-1 < e \leq 2(q-1)$), then the ordered tuple $\mu \in \mathcal{M}_2(n)$ obtained from $\beta$ by 
replacing any entries $\beta_j \leq \beta_{q-1}$ in the range $1 \leq j \leq q-1$ by $\beta_{q-1} -1$ (respectively, replacing any entries $\beta_j \leq \beta_q$ in the range $q \leq j \leq 2(q-1)$ by $\beta_q - 1$) is maximal for $|\cdot|_2$.
\begin{proof}
Let $n$, $\beta$ and $e$ be as in the statement. We have the obvious inequality $|\beta|_1 \geq |\beta|_2$ for ordered tuples $\beta \in \NN^{2(q-1)}$ with equality if and only if the number of even entries equals $q-1$. This deals with the case $e = q-1$.

We give a proof in the case $0 \leq e < q-1$. Explicitly, we are defining 
\[\mu := (\beta_1,\dots,\beta_j, \beta_{q-1} - 1,\dots,\beta_{q-1}-1,\beta_{q},\dots,\beta_{2(q-1)}),\]
where $j \leq e$ is the greatest integer such that $\beta_j \geq \beta_{q-1}-1$.
First, observe that $|\mu|_1 = |\mu|_2$, by definition of $|\cdot|_2$, and we claim that $\mu$ is maximal for $|\cdot|_2$. 
Reduction of any of the entries $\beta_{q},\dots,\beta_{2(q-1)}$ would decrease $|\cdot|_2$, and increasing them would move us out of $\mathcal{M}_2(n)$. Thus we must leave these entries unchanged. We are reduced to analyzing the integer given by 
\[b := q^{\beta_1}+q^{\beta_2}+\cdots+q^{\beta_{q-1}}.\] 
The expansion above gives rise to the base $q$ expansion of $b$ upon collecting the powers of $q$ with equal exponent. From this, we see easily that the greatest integer less than $b$ such that all powers of $q$ have even exponent is
\[m := q^{\beta_1}+\cdots+q^{\beta_j}+(q-1-j)q^{\beta_{q-1} - 1}.\]
The exponents appearing in the equation defining $m$ are exactly the entries of $\mu$. This finishes the case at hand.

In the case that $q-1 < e \leq 2(q-1)$, it appears that we have to worry about entries of $\beta$ that are equal to zero, but this possibility is ruled out by the previous lemma. The remainder of the proof is as in the previous case.

\end{proof}
\end{proposition}

\noindent \textit{Acknowledgements.} The author thanks David Goss and Federico Pellarin who have read and offered comments on more than their fare share of versions of this text. My thanks also go to Daniel Glasscock for comments on a portion of this present text. Finally, I thank the anonymous referee for several useful comments and the encouragement to rewrite the original manuscript, which, I believe, has resulted in a much better paper.

\end{document}